\documentclass[12pt,a4paper]{article}

\def\marginpar#1{}
\def\resizebox#1#2#3{}

\usepackage{bbm}
\usepackage[dvips]{graphicx}
\newcommand{\One}{\mathbbmss{1}}

\usepackage{amsmath}
\usepackage{amsfonts}
\usepackage{amssymb}
\usepackage[mathcal]{euscript}
\usepackage{amsopn}
\usepackage{amsthm}
\usepackage{amstext}

\usepackage{graphics}

\usepackage{smart}

\SectioningParameter{presecnum}{section}=\noexpand\S;
\RestoreLaTeXeqno
\link{equation}{section}

\newcommand {\supplus}{\mathop{{\supset}\llap{\raise
0.5pt\hbox{\normalfont\small+}\hskip 0.5pt}}}

\newcommand {\subplus}{\mathop{{\subset}\llap{\raise
0.5pt\hbox{\normalfont\small+}\hskip 0.5pt}}}

\newcommand {\Cee}    {{\mathbb  C}}

\newcommand {\fg}     {{\mathfrak{g}}}    %
\newcommand {\fgl}    {{\mathfrak{gl}}}  %
\newcommand {\fh}     {{\mathfrak{h}}}

\newcommand {\fn}     {{\mathfrak{n}}}

\newcommand {\fsl}    {{\mathfrak{sl}}}

\def \opname#1#2%
  {\expandafter\newcommand \csname #1\endcsname {{\mathop{#2}\nolimits}}}

\newcommand{\rmname}[1]
  {\expandafter\newcommand \csname #1\endcsname {{\operatorname{#1}}}}

\newcommand{\rmnameii}[2]
  {\expandafter\newcommand \csname #1\endcsname {{\operatorname{#2}}}}

\rmname{act}
\rmname{Ad}
\rmname{Add}
\rmname{ad}
\rmname{Alt}
\rmname{alt}
\rmname{Ann}
\rmname{antidiag}
\rmname{Ber}
\rmname{ber}
\rmname{Br}
\rmname{card}
\rmname{ch}
\rmname{Char}
\rmname{cem}
\rmname{cj}
\rmname{Cliff}
\rmname{cntr}
\rmname{codim}
\rmname{coind}
\rmname{const}
\rmname{col}
\rmname{cork}
\rmname{cpr}
\rmname{diag}
\rmnameii{Div}{div}
\rmname{Def}
\rmname{Der}
\rmname{Dim}
\rmname{End}
\rmname{Even}
\rmname{Ext}
\rmname{gr}
\rmname{Hom}
\rmname{HT}
\rmnameii{Ht}{ht}
\rmname{hwt}
\rmname{Id}
\rmname{id}
\rmname{ind}
\rmname{Ind}
\rmname{Coind}
\rmname{Inf}
\rmname{irr}
\rmname{Le}
\rmname{Lie}
\rmname{lwt}
\rmname{mult}
\rmname{Mor}
\rmname{nm}
\rmname{Ob}
\rmname{Odd}
\rmname{Osc}
\rmname{per}
\rmname{Pic}
\rmname{pr}
\rmname{pro}
\rmname{Prime}
\rmname{Proj}
\rmname{prt}
\rmname{pt}
\rmname{Q}
\rmname{qet}
\rmname{qtr}
\rmname{rd}
\rmname{rk}
\rmname{row}
\rmname{Res}
\rmname{salt}
\rmname{Sch}
\rmname{SBr}
\rmname{scalar}
\rmname{Ser}
\rmname{sign}
\rmname{Smbl}
\rmname{spin}
\rmname{ssym}
\rmname{str}
\rmname{st}
\rmname{sgn}
\rmname{sq}
\rmname{symm}
\rmname{supp}
\rmname{Supp}
\rmname{St}
\rmname{Spec}
\rmname{Spm}
\rmname{tr}
\rmname{vpt}
\rmname{weyl}
\rmname{Weyl}
\rmname{Witt}

\opname{vvol}  {{v\hspace{-0.1ex}o\hspace{-0.02ex}l\/}}

\opname{pnt}  {\text{\normalfont pt}}

\opname{Span} {{Span}}

\opname{slim} {\overline{\lim}}

\opname{Vol}  {{V\hspace{-0.55ex}o\hspace{-0.02ex}l\/}}

\opname{QVol} {{Q\hspace{-0.3ex}V\hspace{-0.55ex}o\hspace{-0.02ex}l\/}}

\opname{PoVol}{{P\hspace{-0.35ex}o\hspace{-0.25ex}V\hspace{-0.55ex}o\hspace{-0.0
2ex}l\/}
}

\opname{BVol} {{B\hspace{-0.2ex}V\hspace{-0.55ex}o\hspace{-0.02ex}l\/}}

\opname{Par}  {{P\hspace{-0.3ex}a\hspace{-0.05ex}r\/}}

\rmname{Mat}
\rmname{Bil}
\rmname{Diff}
\rmname{Ker}
\rmname{Herm}
\rmname{Coker}
\rmname{Conn}
\rmname{Covect}
\rmname{Vect}
\rmname{Int}
\rmname{rank}
\rmname{Pty}

\rmnameii {IM} {Im}
\rmnameii {RE} {Re}

\opname{Aut} {{A\hspace{-0.2ex}u\hspace{-0.1ex}t\/}}
\opname{GL} {{G\hspace{-0.3ex}L}}
\opname{SL} {{S\hspace{-0.3ex}L}}
\opname{Exp} {{E\hspace{-0.2ex}x\hspace{-0.1ex}p\/}}
\opname{GQ} {{G\hspace{-0.2ex}Q}}
\opname{OSp} {{O\hspace{-0.25ex}S\hspace{-0.15ex}p\/}}
\opname{Out} {{O\hspace{-0.25ex}u\hspace{-0.15ex}t\/}}
\opname{Spp} {{S\hspace{-0.2ex}p\/}}
\opname{SpO} {{S\hspace{-0.2ex}p\hspace{-0.02ex}O\/}}
\opname{Pe} {{P\hspace{-0.25ex}e\/}}
\opname{SPe} {{S\hspace{-0.25ex}P\hspace{-0.25ex}e\/}}
\opname{Spin} {{S\hspace{-0.25ex}p\hspace{-0.05ex}i\hspace{-0.1ex}n\/}}
\opname{Iso} {{I\hspace{-0.25ex}s\hspace{-0.1ex}o\/}}
\opname{SSPe} {{S\hspace{-0.25ex}S\hspace{-0.15ex}P\hspace{-0.25ex}e\/}}
\opname{PeU} {{P\hspace{-0.25ex}e\hspace{-0.1ex}U\/}}
\opname{QU} {{Q\hspace{-0.15ex}U\/}}
\opname{U} {{U\/}}

\opname{cGQ} {{\cal G \hspace{-0.2em} Q \/}}
\opname{cSL} {{\cal S \hspace{-0.2em} L \/}}
\opname{cGL} {{\cal G \hspace{-0.2em} L \/}}\opname{cGr} {{\cal G
\hspace{-0.2em} r \/}}
\opname{cOSp} {{\cal O \hspace{-0.2em} S \hspace{-0.3em} \it p\/}}
\opname{cPe} {{\cal P \hspace{-1.5pt} \it e\/}}
\opname{cVect} {{\cal V \hspace{-1.5pt} \it
e\hspace{-0.1ex}c\hspace{-0.1ex}t\/}}
\opname{cVol} {{\cal V \hspace{-1.5pt} \it o\hspace{-0.1ex}l\/}}
\opname{cAut} {{\cal A \hspace{-0.2em} \it u\hspace{-0.1em}t\/}}
\opname{cCovect} {{\cal C \hspace{-1.5pt}
     \it o\hspace{-0.1ex}v\hspace{-0.1ex}e\hspace{-0.1ex}c\hspace{-0.1ex}t\/}}
\opname{CW} {{C\hspace{-0.15ex}W}}

\newcommand {\tto} {\longrightarrow}

\newcommand {\bcdot}   {\mathbin{\hbox{\raise.4ex\hbox{\bf.}}}} 

\newtheorem{Theorem}{Theorem}[section]\toheight2
\newtheorem{Proposition}[Theorem]{Proposition}
\newtheorem{Lemma}[Theorem]{Lemma}

\newtheorem{Corollary}[Theorem]{Corollary}

\theoremstyle{definition}

\theoremstyle{remark}

\begin{document}

\title{Sylvester-'t Hooft generators of $\fsl(n)$ and $\fsl(n|n)$, and relations between them}

\author{Christoph Sachse
\thanks{I am thankful to the International Max Planck Research School for financial support and most
creative environment and to Dimitry Leites for raising the problem.}}

\address{MPIMiS, Inselstr. 22, Leipzig D-04013,
Germany; sachse@mis.mpg.de}



\maketitle

\begin{abstract} Among the simple finite dimensional Lie algebras, only $\fsl(n)$ possesses two
automorphisms of finite order which have no common nonzero
eigenvector with eigenvalue one. It turns out that these
automorphisms are inner and form a pair of generators that allow one
to generate all of $\fsl(n)$ under bracketing. It seems that
Sylvester was the first to mention these generators, but he used
them as generators of the associative algebra of all $n\times n$
matrices $Mat(n)$. These generators
appear in the description of elliptic solutions of the
classical Yang-Baxter equation, orthogonal decompositions of
Lie algebras, 't Hooft's work on confinement operators in QCD, and various other instances.
Here I give an algorithm which both
generates $\fsl(n)$ and explicitly describes a set of defining
relations. For simple (up to center) Lie superalgebras, analogs of
Sylvester generators exist only for $\fsl(n|n)$. The relations for
this case are also computed.
\end{abstract}

\section{Introduction}

Dealing with a given Lie algebra $\fg$ and modules over it,
especially when $q$-quantizing, we need a convenient {\it
presentation} of $\fg$, i.e., a description in terms of generators
and defining relations. Obviously, the basis elements qualify as
generators, but there are too many of them. It is well-known
\cite{GL2} that
\begin{equation}\label{max}
\begin{array}{l}
\text{{\sl For any nilpotent Lie algebra $\fn$, the natural set
of relations is}}\\
\text{{\sl  a basis of $\fn/[\fn, \fn]=H_1(\fn)$; relations
between these generators}}\\
\text{{\sl can be described in terms of the basis of $H_2(\fn)$.
}}
\end{array}
\end{equation}
A simple Lie (super)algebra $\fg$ (finite dimensional, Kac-Moody
or of polynomial vector fields) is conventionally split into the sum
$\fg=\fn_-\oplus\fh\oplus\fn_+$ of two maximal
nilpotent subalgebras $\fn_{\pm}$ (positive and negative) and the
commutative Cartan subalgebra; the corresponding generators are
called {\it Chevalley generators}; the relations between them are
also known, cf. \cite{GL3}, \cite{GLP}. They are numerous ($3n$
generators for a $\rank \,\, n$ algebra and $\sim n^2$ relations), but
these relations are simple and therefore convenient.

For comparison: for the simplest case, $\fgl(n)$, the matrix units
are obvious generators, and the relations between them are simple,
but far too numerous ($n^2$ generators and $\sim n^4$
relations).

Jacobson was, perhaps, the first to observe that every
simple finite dimensional Lie algebra can be generated by just a {\it pair} of
generators, but he did not specify his pairs, so no discussion of relations  was made.
Grozman and Leites \cite{GL1} introduced a
pair of generators associated with the
principal embedding of $\fsl(2)$, and the relations between them
are rather simple (at least, for computers). There are more
generators similar to those Grozman and Leites had chosen, but
experiments performed so far show that the ones Grozman and Leites
considered are most convenient, and are related to various applications
\cite{GL2}, \cite{LS}.

There are, however, certain pairs of generators indigenous only to
the $\fsl$ series, and only over an algebraically closed field, e.g. $\Cee$.
Below, we describe such a pair of generators
for $\fsl(n)$ and their analogs for $\fsl(n|n)$ and give relations between them.

Let $a=\exp{\left(\frac{2i\pi}{n}\right)}$ and define \emph{Sylvester's generators}
(also called {\it clock-and-shift or 't Hooft matrices}) to be
\begin{equation}
D=\mathrm{diag}(1,a,a^2,\ldots,a^{n-1}),\qquad
S=\begin{pmatrix}
0&1 &0&0&0\cr
0&\ddots&\ddots&\ddots&0\cr
0&\ddots&\ddots&\ddots&0\cr
0&\ddots&\ddots&\ddots&1\cr
1&0&0&0&0
\end{pmatrix}
\label{sylv}
\end{equation}

Zachos \cite{Z1} points out that \\[1em]
\lq\lq {\sl apparently, Sylvester \cite{S} was the first to study
these (\ref{sylv}) generators\footnote{More precisely, Sylvester
used them as generators of an associative algebra, where they
yield the algebra of $n\times n$ matrices $\Mat(n)$. Having
replaced the dot product by the bracket we endow the space of
$\Mat(n)$ with the structure of the Lie algebra $\fgl(n)$; having
introduced parity in $\Mat(n)$ by attributing parity to each basis
vector (and hence to each row and column) and replacing the dot
product by the superbracket we endow the superspace of $\Mat(n;
\Par)$, where $\Par$ is an ordered collection of parities, with
the structure of the Lie superalgebra $\fgl(\Par)$. As generators
of a Lie algebra or Lie superalgebra, Sylvester's generators can
only generate $\fsl(n)$ and $\fsl(n|n)$, not $\fgl(n)$.} of
$\fsl(n)$; he worked them out for $\fsl(3)$ first, and called them
\lq\lq nonions" (after quaternions), and then generalized to
$\fsl(n)$.

They became popular in the 30s in the context of QM-around-the circle,
i.e., on a discrete periodic lattice of $N$ points, see \cite{W}.
That effort has continued to date, with the work of Schwinger,
Santhanam, Tolar, Floratos, and others.

They also became popular among high-energy theorists, with the
work of 't~Hooft \cite{tH}, on order-disorder confinement
operators in QCD, so that many in my end of the woods intriguingly
call them \lq\lq 't Hooft matrices".

I have been using them every few years, starting from \cite{FFZ}
to identify cases of a Sine-algebra we found at that time
with $\fsl(N)$, and also with the Moyal Bracket algebra
\cite{Moy} on a toroidal phase space; and hence take the $N\tto \infty$
limit
to get Poisson Brackets more directly than in Hoppe's first derivation
\cite{Ho} on a spherical phase space.

Our latest use of them was in our recent diversion, \cite{FZ}, on
ring-indexed Lie algebras. They are apparently the most systematic
basis for dealing with all $\fsl(N)$s on an equal footing and taking
naive $N\tto \infty$ limits.}"\\[1em]

For the passage from the notation of Zachos et al.
to ours, observe that, e.g. in \cite{FFZ}, the authors generate $\fgl(n)$ from
Sylvester's generators $D,S$ (\ref{sylv}) in the form
\[
J_{(m_1,m_2)}=a^{m_1m_2/2}D^{m_1}S^{m_2}
\]
which are $n^2$ independent matrices labelled by two integers $0\leq m_1,m_2<n$. Under the
bracket, the identity matrix $J_{(0,0)}$ spans the center. So dividing it out leaves
$\fsl(n)$ with the bracket
\[
[J_{(m_1,m_2)},J_{(k_1,k_2)}]=-2i\sin\left( \frac{2\pi}{n}(m_1k_2-m_2k_1) \right) J_{(m_1+k_1,m_2+k_2)}
\]

Another important
application of Sylvester's generators is the classical Yang-Baxter equation for a function taking
values in a simple Lie algebra $\fg$. It turns out \cite{BD1, BD2} that for this equation to
have elliptic solutions, $\fg$ has to possess two automorphisms
of finite order which have no common nonzero eigenvector with
eigenvalue 1. Sylvester's generators are such automorphisms for $\fg=\fsl(n)$;
in fact, \cite{BD1, BD2} prove that any $\fg$ possessing such automorphisms must
be isomorphic to $\fsl(n)$, and the elliptic solutions can
be characterised by the images of Sylvester's generators (\ref{sylv}) under this isomorphism.
Also, they play a vital role in the study of orthogonal
decompositions of Lie algebras \cite{KKU, KT, FOS}.

Finally, a more applied subject on which these generators have been used is
hydrodynamics and the statistical theory of turbulent fluids and gases, in particular,
the study of lattice models of inviscid fluids (Euler fluids), see, e.g., \cite{MWC},\cite{MW},\cite{Ze}.

The aim of this paper is to give an algorithm that generates $\fsl(n)$ and $\fsl(n|n)$ from
Sylvester's generators and which also produces a presentation for them. This presentation contains
redundancies, but might be of interest for practical problems since it allows quick and easy
computations in the adjoint representation. The main statements are the following ones.

\begin{Theorem}
Fix an integer $n\geq 2$. Then the matrices (\ref{sylv}) are generators for $\fsl(n)$:
\begin{equation}
\mathfrak{sl}(n)=Span(D,S,T_m^k\mid 1\leq k,m \leq n, \,\,\mathrm{and}\,\,k\neq n\,\,%
\mathrm{for}\,\,m=1,n,\,\,\mathrm{and}\,\,k\neq 1\,\,\mathrm{for}\,\,m=n),
\label{spansl}
\end{equation}
where for $1\leq k,m\leq n$, we set
\begin{eqnarray*}
T_m^k &=& (\ad\,D)^{k-1}((\ad\,S)^{m-1}((\ad\,D(S))),\\
T_n^k &=& \ad\,S(T_{n-1}^k).
\end{eqnarray*}
A defining set of relations for generators (\ref{sylv}) can be
obtained in the following way. The relations
\begin{eqnarray}
\label{rel1}
(\ad D)^n(S) &=& (1-a)^n S,\\
\label{rel2}
(\ad D)^{n}((\ad S)^{m-1}(\ad D(S))) &=& (1-a)^m(1-a^m)^n(-1)^{m+1}(\ad S)^{m-1}(\ad D(S))\\
\label{rel3}
\ad S(T_{n-1}^1) &=& (1-a)^n(-1)^{n}D,\\
\label{rel4}
\ad S(T_{n-1}^n) &=& 0,\\
\label{rel5}
\ad S(T_n^k) &=& (-1)^n(1-a^k)^2(1-a^{n-1})^{k-1}(1-a)^{n-k}T_1^k,\\
\label{rel6}
\ad D(T_n^k) &=& 0
\end{eqnarray}
prohibit generation of elements of order higher than $n$ in both
$D$ and $S$. Besides them, for each $T_m^k$ with $2\leq m\leq
n-1$, except for $T_2^2$, $m-1$ relations have to hold, which can
be written as
\[
(\ad\,S)^{s_1}((\ad\,D)^{k-1}((\ad\,S)^{s_2}(T_1^1)))=\left(\frac{1-a^k}{1-a}\right)^{s_1}%
\left(\frac{1-a^{s_1}}{1-a^{m-1}}\right)^{k-1}T^k_m,
\]
where $s_1+s_2=m-1$ and $s_1=1,2, \ldots, m-1$. \label{main1}
\end{Theorem}

\begin{Theorem}
Considered as $2n\times 2n$ supermatrices on a superspace with an
alternating format (even, odd, even, odd, ...), (\ref{sylv}) are
generators for $\fsl(n|n)$:
\[
\mathfrak{sl}(n|n)=Span(D,S,T_m^k\mid 1\leq k,m \leq 2n, \,\,\mathrm{and}\,\,k\neq 2n\,\,%
\mathrm{for}\,\,m=1\,\,\mathrm{and}\,\,k\neq 1\,\,\mathrm{for}\,\,m=2n)
\]
with the same definition of $T_m^k$ as in Thm. (\ref{main1}).
A defining set of relations in this case are
(\ref{srel1})-(\ref{srel5}) and $m-1$ relations for each $T_m^k$ with $2\leq m\leq 2n-1$.
These can be written as
\[
(\ad\,S)^{s_1}(T^k_{s_2+1})=\begin{cases}
(-1)^{s_1}(a^{2k}-1)^{\frac{s_1-1}{2}}\frac{(a^k+1)}{1-a}\left(\frac{1-a^{s_1}}{1-a^{m-1}}\right)^{k-1}%
\left(-\frac{1-a}{1+a}\right)^{\frac{s_1-1}{2}}\tilde{T}_m^k & \text{for $s_1$ odd and}\\
& \text{$s_2$ even,}\\
(-1)^{s_1}(a^{2k}-1)^{\frac{s_1-1}{2}}\frac{(a^k-1)}{1-a}\left(\frac{1-a^{s_1}}{1-a^{m-1}}\right)^{k-1}%
\left(-\frac{1-a}{1+a}\right)^{\frac{s_1-1}{2}}\tilde{T}_m^k & \text{for $s_1,s_2$ odd, or}\\
(-1)^{s_1}\frac{(a^{2k}-1)^{s_1/2}}{1-a}\left(\frac{1-a^{s_1}}{1-a^{m-1}}\right)^{k-1}%
\left(-\frac{1-a}{1+a}\right)^{s_1/2}\tilde{T}_m^k & \text{for $s_1,s_2$ even,}\\
& \!\!\!\!\!\text{or $s_1$ even, $s_2$ odd,}
\end{cases}
\]
where again $s_1+s_2=m-1$ and $s_1=1,2,\ldots,m-1$.
\label{main2}
\end{Theorem}

In the following we show why the set of relations indicated in
Thms. (\ref{main1}), (\ref{main2}) is a defining set. This will
then automatically also deliver an upper bound on the number of
independent relations for Sylvesters's generators.

\section{Relations between Sylvester's generators for $\fsl(n)$}

Setting
\[
T_1^k:=(\ad D)^k(S)\qquad\mathrm{for}\quad k=1,\ldots,n-1
\]
we obtain the matrices
\[
T_1^k=(1-a)^{k}\left(\begin{array}{ccccc}
0 & 1 & 0 & \ldots & 0\\
0 & 0 & a^k & \ldots & 0\\
\vdots\\
0 & 0 & \ldots & \ldots & a^{k(n-2)}\\
a^{k(n-1)} & 0 & \ldots & 0 & 0
\end{array}\right)
\]
which are, clearly, all linearly independent. For $k=n$, we get
the relation (\ref{rel1}). Proceeding likewise, we generate a
basis for $\mathfrak{sl}(n)$. We set
\[
T_m^k=(\ad D)^{k-1}((\ad S)^{m-1}(T_1^1))=(\ad D)^{k-1}((\ad\,S)^{m-1}((\ad\,D(S)))
\]
where $m=1,\ldots,n-1$ and $k=1,\ldots,n$. In matrix form,
\[
T_m^k=(1-a)^m(1-a^m)^{k-1}(-1)^{m+1}\left(\begin{array}{cccccccc}
0 & \ldots & 0 & 1 & 0 & \ldots & \ldots & 0\\
0 & \ldots & 0 & 0 & a^k & 0 & \ldots & 0\\
\vdots & \vdots &&&&&  & \vdots\\
0 & \ldots &&&&&  & a^{k(n-m-1)}\\
a^{k(n-m)} & 0 & \ldots &&&& \ldots & 0\\
\vdots &&&&&&& \vdots\\
0 & \ldots & a^{k(n-1)} & 0 & \ldots && \ldots & 0
\end{array}\right)
\]
These are all the non-diagonal matrices needed for a basis of
$\mathfrak{sl}(n)$. Their linear independence is easily checked.
We also immediately read off the relation (\ref{rel2}) for $k=n+1$.

It remains to generate $n-2$ diagonal matrices, which we do as
follows:
\[
T_n^k=\ad S(T_{n-1}^k)=\ad S((\ad D)^{k-1}((\ad S)^{n-2}(\ad
D(S)))),
\]
where $k=2,\ldots,n-1$, and we obtain the relations (\ref{rel3})-(\ref{rel6}).

Since we have obtained $n^2-1$ linearly independent matrices, we
have found a basis for $\mathfrak{sl}(n)$, see (\ref{spansl}). One
might, however, wish to describe $\mathfrak{sl}(n)$ as the
quotient of the free Lie algebra generated by the two Sylvester
generators modulo certain defining relations. Since we know that
the matrices $T_m^k$ span $\mathfrak{sl}(n)$, we know that any
commutator of them must yield a relation. The relations stated
above are merely those ones that are first encountered when we
proceed through our chosen algorithm for the generation of the
basis of $\mathfrak{sl}(n)$.  To find out the number and an
explicit realization of the minimal defining relations turns out
to be quite a tough job, despite the seeming simplicity of the
problem. P.~ Grozman was able to find those minimal relations for
$n=2,3,4$:
\begin{equation}\label{grrel}
\renewcommand{\arraystretch}{1.3}
\begin{array}{l}
\underline{n=2}: \;(\ad\,S)^2(D)=4D,\quad (\ad\,D)^2S=4S;\\
\underline{n=3}: \;(\ad\,S)^3(D)=-3(a-a^2)D,\quad (\ad\,D)^3(S)=3(a-a^2)S;\quad
[T^1_2, T^2_1]=0;\\
\underline{n=4}: \; (\ad\,S)^4(D)=-4D,\quad (\ad\,D)^4(S)=-4S,\quad [T_1^1, T^1_2]=[D,T^1_3],\\
\quad [T_1^1,[T^1_2, T^1_3]]=-4T^2_1,\quad [T_1^3, T_3^1]=0,\quad%
2[T^1_2, T^3_1]=[T^2_1, [D,T^1_2]],\\
\quad 2[T^2_1, T^1_3]=[T^1_2, [S, T_1^2]], T^2_1]=[S,T^3_1]\quad
[T^2_1, T^3_1]=4T^1_2.
\end{array}
\end{equation}
with the help of {\it Mathematica} and his {\bf SuperLie} package \cite{Gr}, but did not succeed to
deduce from (\ref{grrel}) a
general formula. On the other hand, neither the number nor an
explicit form of the \emph{minimal} set of relations is of great practical
importance when working with these generators. Rather one would
like to have, e.g., formulae that describe the action of arbitrary
products of  the elements of $\mathfrak{sl}(n)$
in the adjoint representation. Such formulae will be given below and,
additionally, a set of relations offered which
contains redundancies, but which allows immediate reduction of an arbitrary expression of
the form (with the $X_i$ and $Y$ being arbitrary elements of $\fsl(n)$)
\[
\ad\,X_1(\ad\,X_2(\ldots(\ad\,X_q(Y))\ldots))
\]
to a linear combination of the basis elements produced by our algorithm.

By explicit calculation one first verifies that
\[
[T_m^k,T_{m'}^{k'}]=\frac{(1-a^m)^{k-1}(1-a^{m'})^{k'-1}}{(1-a^{m+m'})^{k+k'-1}}(a^{k'm'}-1)T^{k+k'}_{m+m'}
\]
for any of the $T_m^k,T_{m'}^{k'}$ defined above. Hereafter, $k+k'$ and $m+m'$
have to be understood $\mathrm{mod}\,\, n$. This directly shows the following statement.

\begin{Lemma} The result of the application of an arbitrary product of elements in the adjoint
representation to a $T_m^k$ depends \emph{up to a factor} only on the number of $S$'s and $D$'s
contained in these operators. That is,
\begin{equation}
\ad \,X_1(\ad\, X_2(\cdots\ad\, X_N (T_m^k))\cdots) = C(k,k',m,m',n)T_{m+m'}^{k+k'}
\label{lem1_stat}
\end{equation}
where $m_i$ and $k_i$ is the number of $S$'s and $D$'s, respectively, contained in $X_i$,
\[
m'=\sum_i m_i,\qquad k'=\sum_i k_i,
\]
and $C(k,k',m,m',n)$ is a constant depending on all the indices.
\label{lem1}
\end{Lemma}

Therefore we conclude that it suffices to check only relations between elements which are at most of
degree $n$ in both $S$'s and $D$'s. If we know all relations of this type, then any relation of a
higher degree will follow from these and (\ref{rel1})-(\ref{rel6}).

In order to find these relations, it is most convenient to visualise the generated basis as
a grid of points.
\begin{figure}
\begin{center}
\includegraphics{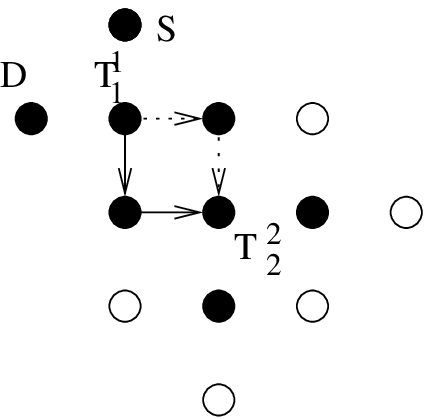}
\caption{The two paths that generate $T_2^2$ in $\mathfrak{sl}(3)$}
\label{fig2p}
\end{center}
\end{figure}
Fig. \ref{fig2p} shows the basis of $\mathfrak{sl}(3)$, starting from $D$ and $S$ in the upper
left corner. Below them is $T_1^1=[D,S]$. The other solid points are those that we generate with our algorithm by
going only horizontally on each level, and vertically only along the left edge. The white points
are those which are ruled out by the relations (\ref{rel1})-(\ref{rel6}), i.e., they do not represent
basis elements of $\mathfrak{sl}(3)$. Now, an arbitrary product of $r$-many $(\ad D)$'s and
$s$-many $(\ad S)$'s applied to
$T_1^1$ corresponds to a path on the grid starting at $T_1^1$ and reaching $T_{s+1}^{r+1}$,
but one which will in general
only produce a matrix proportional to $T_{s+1}^{r+1}$, with a factor $\neq 1$. A horizontal step of the
path
describes the action of $\ad\,D$, a vertical one the action of $\ad\,S$. In the picture, the
solid line shows the way our algorithm went to generate $T_2^2$, while the dotted lines show the
alternative path, i.e.,
\begin{eqnarray}
\textrm{solid line}\qquad & \Leftrightarrow & \quad\ad D(\ad S(\ad D(S))) \label{path1}\\
\textrm{dotted line}\qquad & \Leftrightarrow & \quad\ad S((\ad D)^2(S)) \label{path2}
\end{eqnarray}
It is clear that any expression we have to examine can be expressed as a path from $T_1^1$ to some
admissible $T_m^k$ which only moves right and downwards (compare to Fig. 1). In general, there are
\[
\left(\begin{array}{c}
m+k-2\\
k
\end{array}\right)\qquad\textrm{paths from}\,\, T_1^1\,\,\mathrm{to}\,\,T_m^k
\]
However, we can rule out some of these. The algorithm always uses paths which run through all
vertical steps first, then through all horizontal ones (called the algorithm path in what follows).
A relation is obtained by running through
any different path and comparing the result to what the algorithm path would have produced at this vertex.
\begin{Proposition} For a path ending at $T_s^r$ to yield an independent relation, it has to
\begin{itemize}
\item end with a vertical step if $s<n$,
\item end with a
horizontal step if $s=n$.
\end{itemize}
\label{prop1}
\end{Proposition}
\begin{proof}
Look at the $s<n$ case first.
We know that at the vertex where the last vertical step ends, we will have produced a matrix
proportional to the one that the algorithm path would have produced there (cf. Lemma \ref{lem1}).
Thus, at this vertex we obtain a relation.
But if it is followed by horizontal steps, these will then trivially also yield matrices proportional
to those that the algorithm would have produced. Thus, the relations we can read off at these vertices
are generated from the one obtained at the end of the last vertical step.

An analogous argument holds for $s=n$, except that at the last step of the algorithm there is a vertical
step, so a path producing an independent relation cannot have a vertical step at its end.
\end{proof}

\begin{Corollary}
Apart from those vertical steps which lie on the left edge of the grid,
a path that leads to $T_s^r$ and yields an independent relation for $s<n$ must contain all other
vertical steps at its end. For $s=n$, the only path yielding a nontrivial relation is the algorithm
path to $T_n^{r-1}$ followed by a horizontal step.
\label{cor1}
\end{Corollary}
\begin{proof}
As a counterexample for the $s<n$ case, consider Fig. \ref{dstep}.\\
\begin{figure}[!ht]
\begin{center}
\includegraphics{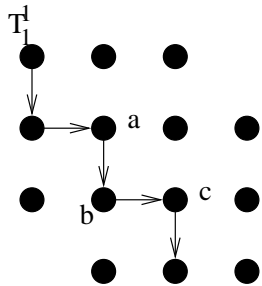}
\caption{Example of a path ruled out by Corollary \ref{cor1}}
\label{dstep}
\end{center}
\end{figure}
\newline
Up to vertex $a$, it follows the algorithm path, then going to $b$ will yield a relation. But
proceeding further horizontally after $b$ yields only dependent relations, as seen before.

In the $s=n$ case, we have seen in Lemma \ref{lem1} that the last step of a path yielding a relation
must be horizontal. Since going a horizontal step in the $n$-th row always gives zero
(cf. (\ref{rel6})), a
nontrivial path can only have exactly one horizontal piece at its end. So the second last step is
always the last step of the algorithm to $T_n^{r-1}$, and therefore any other path leading to
$T_n^{r-1}$ followed by a horizontal step would trivially yield a result proportional to what the
algorithm path followed by the horizontal step gives. The relations so obtained are precisely those of
(\ref{rel6})
\end{proof}

This reduces the number of possibly independent relations considerably: for any vertex $T_m^k$ with $m<n$,
there can now
be at most $m-1$ independent relations, which result from the paths leading there and having
between zero and $m-1$ vertical steps at their ends. For $T_n^k$, there can only be one
relation.  Among the relations thus obtained,
there will still be redundancies, which are not obvious at first glance. To reveal them, one has to
apply the Jacobi identity and other relations one has already obtained. As an example, look at
$T_2^2$ in $\mathfrak{sl}(n)$ for $n\geq 3$. Two paths lead there, described in (\ref{path1})
and (\ref{path2}). Since they are both admissible in the sense of Proposition \ref{prop1}, one might
think that we obtain a relation here between $T_2^2$ generated by the algorithm and
the result of another path. However,
\[
\ad D(\ad S(\ad D(S)))=\ad S((\ad D)^2(S))+\ad(\ad D\,(S))(\ad D(S))
\]
due to the Jacobi identity and the last term is of the form $\ad z(z)\equiv 0$.
Therefore, the two paths \emph{trivially} yield the same result, and we obtain no relation here.
We will show now that there are no other interdependencies of this sort
except the above one for $T_2^2$.\\
\begin{Lemma}
It is impossible to trivially identify the result of two paths to a given $T_s^r$
by rearranging them using the Jacobi identity, except for the case $r=s=2$, where we have
\[
\ad D(\ad S(\ad D(S)))=\ad S((\ad D)^2(S))
\]
\label{comml}
\end{Lemma}
\begin{proof}
Any admissible path in the sense of Lemma \ref{lem1} and its corollary is of the form
\begin{equation}
T_s^r=(\ad\,S)^{s_2}((\ad\,D)^{r-1}((\ad\,S)^{s_1}(T_1^1)))
\label{path}
\end{equation}
where $s_1+s_2=s-1$. For $s_2=0$, we obtain the algorithm path.
In order to show that two paths give the same result, we want to apply the Jacobi identity
\[
\ad(\ad\,x\,\,(y))z=\ad\,x\,(\ad\,y\,\,(z))-\ad\,y\,(\ad\,x\,\,(z))
\]
in such a way that the left hand side of it becomes zero, i.e. is of the form
$\ad\,z(z)$. This would rule out one of the relations these paths produce. We see
immediately that for this to happen for adjoint operators $x,y,z$, the element $z$
would have to contain as
many $D$'s and $S$'s as $x$ and $y$ together. Looking at (\ref{path}), which we would like to
identify with $\ad\,x\,(\ad\,y\,\,(z))$, this implies $r=s$. We have to split (\ref{path}) in two
equally long subpaths, the head (including $T_1^1$) being $z$ and the tail being $\ad\,x\,(\ad\,(y))$ and
each containing $\frac{s}{2}$-many $D$'s and $S$'s, implying $s$ must be even.

For the case $r=s=2$, we find that $z=T_1^1$, $x=\ad\,S$ and $y=\ad\,D$ meet these requirements.

Let now $r=s=2n$, $n>1$ and let $x,y,z$ satisfy the above conditions.
Then $z$ represents the path of the algorithm to
$T_{s/2}^{s/2}$ and $\ad\,x(\ad\,y)$ is of the form $(\ad\,S)^{s/2}((\ad\,D)^{s/2})$. But we
see that it is impossible then to find $x,y$ such that $\ad\,y((\ad\,x)(z))$ would again be an
admissible path.
\end{proof}

It is important to note that this still does not exclude all possible dependencies between the
relations that various admissible paths yield. By clever rearrangement, it might still be possible
to bring a bracket of two elements into a form which, when expanded into paths, yields only a few admissible
paths and several others which run over already excluded pieces. We could find no way to rule out
all such possibilities. This seems only possible with the help of computers.
But, as stated above, the minimal number might not be of practical interest. The preceding discussion
still gives us an upper bound on the number of relations.\\
\begin{Theorem}
The number $R(n)$ of independent relations between Sylvester's generators is
bounded from above by:
\[
R(n)\leq\begin{cases}
2 & \mathrm{for}\,\,n=2,\\
n^2-3 & \mathrm{for}\,\,n\geq 3.
\end{cases}
\]
\end{Theorem}
\begin{proof}
\underline{$n=2$:}
See Fig. \ref{n2fig} for the diagram.
\begin{figure}[!ht]
\begin{center}
\includegraphics{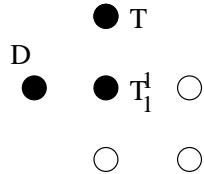}
\caption{The grid of basis elements for $\mathfrak{sl}(2)$}
\label{n2fig}
\end{center}
\end{figure}
\newline
The white dots are ruled out by the relations stated in the beginning, however the dot in the
lower right corner is not independent here. Thus, the only relations are
$$
(\ad D)^2(S)=4S,\qquad
(\ad S)^2(D)=4D.
$$
which was also Grozman's result (\ref{grrel}).\\
\underline{$n\geq 3$:}
As an example for the generic case, look at the $n=5$ grid (Fig. \ref{n5fig}).
\begin{figure}[!ht]
\begin{center}
\includegraphics{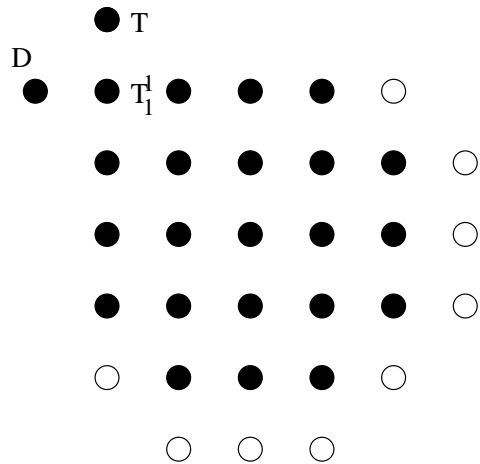}
\caption{The grid of basis elements for $\mathfrak{sl}(5)$}
\label{n5fig}
\end{center}
\end{figure}
\newline
We get here the following relations:
\begin{itemize}
\item 2 relations for the $n$-th powers of $\ad D$ and $\ad S$,
\item $(n-2)$ relations that limit the application of $\ad\,D$ (the rightmost white dots),
\item $(n-2)$ relations that limit the application of $\ad\,S$ (lowermost white dots),
\item 1 relation corresponding to relation (\ref{rel4}) (white dot in the lower right corner),
\item $(n-3)$ relations for the vertical paths from the first to the second row,
\item $(n-3)(n-1)$ relations for the vertical paths between the second and third row, third and
fourth row and so on down to the $(n-1)$st row,
\item $(n-2)$ relations for the horizontal paths in the $n$-th row.
\end{itemize}
This makes a total of $n^2-3$ relations. Lemma \ref{lem1} and its corollary exclude the possibility
that one of them is obtained by another by application of $\ad\,D$ or $\ad\,S$. Lemma \ref{comml}
shows that none of them is a consequence of another via a rearrangement using the Jacobi identity.
\end{proof}

We see that even for $n=3$, the bound overestimates the exact number of relations. However, the
number of relations found in the above manner is only of order $\sim n^2$, which can be expected to lie
pretty close to the true behaviour of $R(n)$ so that the relative error will decrease for
growing $n$. But the main advantage of our method is that it
explicitly produces a \emph{presentation} (albeit a redundant one): all relations can be directly read
off from the grid representation of the basis of $\fsl(n)$.

\section{Relations between Sylvester's generators for $\fsl(n|n)$}

Sylvester's generators can as well be used to generate a basis of $\fsl(n|n)$, and only for this
simple (up to a nontrivial center) finite dimensional Lie superalgebra, see \cite{LSe}. It is
most convenient to choose an alternating format for the superspace in which we express
the supermatrices, i.e., if $(e_1,e_2,\ldots,e_n)$ is a basis of this vector space, let the $e_{2k+1}$
be
odd vectors and the $e_{2k}$ be even ones for all $k$. This format has the advantage that we can use
the same matrices $D,S$ as above as Sylvester's generators, where now $D$ is an even supermatrix and
$S$ an odd one. The result obtained below remains valid in any format, but looks nicest in the chosen one.
To be able to compare the matrices obtained for the $\fsl(n)$ and $\fsl(n|n)$ cases, we put a twiddle
on the supermatrices: $\tilde{D},\tilde{S}$.

As above, set
\[
\tilde{T}_1^1=[\tilde{D},\tilde{S}]
\]
which is now an odd supermatrix, but with the same entries as in the $\fsl(n)$ case. Likewise,
\[
\tilde{T}_1^k=(\ad\,\tilde{D})^{k-1}(\tilde{T}_1^1)
\]
are all odd supermatrices, but look the same as in the $\fsl(n)$ case, and we find the analogue of relation
(\ref{rel1}) to be
\begin{equation}
(\ad\, \tilde{D})^{2n}(\tilde{S})=(1-a)^{2n}\tilde{S}.
\label{srel1}
\end{equation}
We follow the same algorithm as in the $\fsl(n)$ case, now using the superbracket: set
\[
\tilde{T}_2^1=[\tilde{S},[\tilde{D},\tilde{S}]]=[\tilde{S},\tilde{T}_1^1]
\]
which is the same matrix as in the $\fsl(n)$ case, except for the prefactor, which is now $(1-a^2)$
instead of $-(1-a)^2$. In general, for $k=1,\ldots,2n;\,\,m=1,\ldots,2n-2$,  we have
\begin{eqnarray*}
(\ad\,\tilde{D})^{k-1}([\tilde{D},\tilde{S}]) &=& \tilde{T}^k_1 = T^k_1\qquad\mathrm{for}\quad %
k=1,\ldots,2n-1\\
(\ad\,\tilde{D})^{k-1}((\ad\,S)(\tilde{T}^1_m)) &=& \tilde{T}^k_{m+1} = %
\left\{\begin{array}{ll}
-\left(-\frac{1+a}{1-a}\right)^{\frac{m+1}{2}}T^k_{m+1} & \mathrm{for}\,\,m\,\,\mathrm{odd}\\
-\left(-\frac{1+a}{1-a}\right)^{\frac{m}{2}}T^k_{m+1} & \mathrm{for}\,\,m\,\,\mathrm{even}
\end{array}\right.\\
\end{eqnarray*}
so $\tilde{T}^k_m$ is proportional to $T^k_m$. The $(1+a)$-factors stem from the
application of anticommutators. One obtains the analogue of the
relations (\ref{rel2}) for $k=2n+1$ and $2\leq m \leq 2n-1$:
\begin{equation}
(\ad \tilde{D})^{2n}((\ad \tilde{S})^{m-1}((\ad \tilde{D})(\tilde{S})))=%
\left(-\frac{1+a}{1-a}\right)^{\frac{m(+1)}{2}}(1-a^m)^{2n}(-1)^{m}%
(\ad \tilde{S})^{m-1}(\ad\,\tilde{D}(\tilde{S})).
\label{srel2}
\end{equation}
For the diagonal basis elements, we set
\[
\tilde{T}_n^{k}=\ad\,\tilde{S}(\tilde{T}_{n-1}^k)\qquad\mathrm{for}\quad 2\leq k\leq n
\]
and obtain the following relations:
\begin{eqnarray}
\label{srel3}
\ad\, \tilde{S}(\tilde{T}_{n-1}^1) &=& \left(-\frac{1+a}{1-a}\right)^n(-1)^{n+1} \tilde{D},\\
\label{srel4}
\ad\,\tilde{S}(\tilde{T}_n^k) &=& %
\left(-\frac{1+a}{1-a}\right)^{n-1}(1-a)(a^{2k}-1)(1-a^{2n-1})^{k-1}(-1)^{n}\tilde{T}_1^k.
\end{eqnarray}
Note that $\tilde{T}^{2n}_{2n}$ is not zero here, but is proportional to the identity matrix. On any
$(n|n)$-dimensional superspace, the identity matrix is supertraceless, and therefore an element of $\fsl(n|n)$.
Thus, no relation corresponds to (\ref{rel4}) in the super case.

We have to add one more relation, which did not exist in the non-super case: the supercommutator of
$\tilde{S}$ with itself:
\begin{equation}
[\tilde{S},\tilde{S}]=\frac{1}{(1+a)(1-a)^{2n-1}}\tilde{T}^{2n}_2\qquad\text{for $n>1$.}
\label{srel5}
\end{equation}
For $n=1$, this is not a relation, but really generates a new element, see Thm.
(\ref{superthm}).

Thinking of the set of basis elements again as a grid of points, we see that we have found relations of the
same sort as in the $\fsl(n)$ case, with one exception: there is one more element, the one proportional
to the identity matrix, represented by the rightmost dot in the last row.

One can again verify by explicit calculation that $[\tilde{T}^k_m,\tilde{T}^{k'}_{m'}]$ is proportional
to $\tilde{T}^{k+k'}_{m+m'}$. This extends the validity of Lemma \ref{lem1} to the super case.
To find a bound for the number of relations again reduces to checking all paths from $\tilde{T}_1^1$
to the other $\tilde{T}_m^k$'s. This is done in the same way as before, it is clear that our algorithm
proceeds on the same paths as in the $\fsl(n)$ case and that Prop. \ref{prop1} and Cor. \ref{cor1} also
apply in the super case.

Also Lemma \ref{comml} generalises to the super case, now using the super Jacobi identity. But here we have
to be careful about a specialty of the super case: supercommutators of elements with themselves do not
necessarily vanish. Consider, for example, $\tilde{T}_2^2$:
\begin{equation}
\ad\,\tilde{D}(\ad\,\tilde{S}(\ad\,\tilde{D}\,(\tilde{S})))=%
-\ad(\ad\,\tilde{S}(\tilde{D}))(\ad\,\tilde{D}(\tilde{S}))-\ad\,\tilde{S}((\ad\,\tilde{D})^2(\tilde{S})).
\end{equation}
Here, the first term on the right hand side does not vanish. Therefore the relation between the two paths
to $\tilde{T}_2^2$ that we ruled out as being trivial in the $\fsl(n)$ case is nontrivial in the super
case. Except for this fact, Lemma \ref{comml} remains valid.

\begin{Theorem}
For $\fsl(n|n)$, the number $R(n)$ of independent relations between Sylvester's generators is bounded by
\begin{equation}
R(n)\leq\begin{cases}
4 & \text{for $n=1$,}\\
(2n)^2-1 & \text{for $n>1$.}
\end{cases}
\end{equation}
\label{superthm}
\end{Theorem}
\begin{proof}
The $n=1$ case differs from the $\fsl(2)$ case because of the relation
\[
[\tilde{S},\tilde{S}]=2\cdot\One,
\]
where $\One$ is the identity matrix.
\begin{figure}[!ht]
\begin{center}
\includegraphics{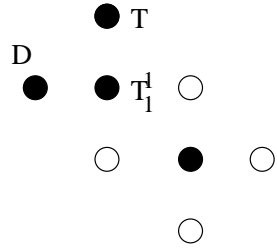}
\caption{The grid of basis elements of $\fsl(1|1)$}
\label{fig11}
\end{center}
\end{figure}
The basis elements of $\fsl(1|1)$ can be represented by the grid in
Fig \ref{fig11}.
There are four relations:
\begin{eqnarray}
[\tilde{D},[\tilde{D},\tilde{S}]] &=& 4\tilde{S}\\
{}[\tilde{S},[\tilde{D},\tilde{S}]] &=& 0\\
{}[\tilde{D},[\tilde{S},\tilde{S}]] &=& 0\\
{}[\tilde{S},[\tilde{S},\tilde{S}]] &=& 0
\end{eqnarray}

\begin{figure}[!ht]
\begin{center}
\includegraphics{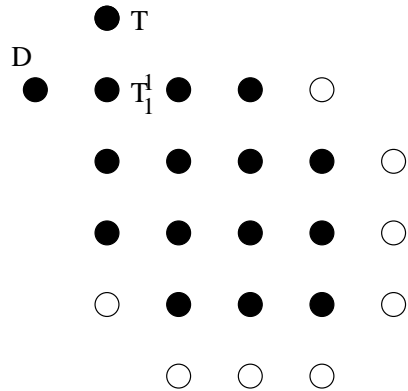}
\caption{The grid of basis elements of $\fsl(2|2)$}
\label{fig22}
\end{center}
\end{figure}
For $n>1$, the grid looks like in Fig. \ref{fig22}.
Note that now there is one more black dot in the lower right corner which we generate from the dot above
it. This provides one more relation. Another additional relation is obtained from the two paths to
$\tilde{T}_2^2$, which are now independent. Apart from this, the situation is identical to the non-super
case.
\end{proof}

\end{document}